\documentclass[12pt]{amsart}

\usepackage[cp1250]{inputenc}

\usepackage{amsmath}
\usepackage{amssymb}
\usepackage{amsthm}
\usepackage{setspace}

\usepackage[varg]{pxfonts}

\numberwithin{equation}{section}
\newtheorem{thm}{Theorem}[section]
\newtheorem{lem}[thm]{Lemma}
\newtheorem{cor}[thm]{Corollary}

\newcommand{\set}[1]{\left\{{#1}\right\}}
\newcommand{\abs}[1]{\left\vert{#1}\right\vert}

\newcommand{\snorm}[1]{\Vert{#1}\Vert}
\newcommand{\C}{\mathbb C}
\newcommand{\R}{\mathbb R}

\renewcommand{\H}{\mathcal H}

\newcommand{\U}{\mathcal U}

\newcommand{\Ci}{C^\infty}
\newcommand{\la}{\langle}

\newcommand{\ra}{\rangle}
\newcommand{\seq}{\bar}
\newcommand{\Cl}[1]{\overline{#1}}
\newcommand{\res}[1]{|_{#1}}
\newcommand{\X}{\mathfrak{X}}
\newcommand{\Ld}{\mathcal{L}}
\newcommand{\ip}{\mathbin{\lrcorner}}
\renewcommand{\Re}{\mathop{\mathrm{Re}}}
\newcommand{\Fl}{\mathrm{Fl}}
\newcommand{\conj}[1]{\overline{#1}}

\newcommand{\pmt}{\Pi_\mu^\theta}

\newcommand{\pmz}{\Pi_\mu^0}

\newcommand{\supp}{\mathop{\mathrm{supp}}}
\newcommand{\Sympl}{\mathop{\mathrm{Sympl}}}
\newcommand{\Diff}{\mathop{\mathrm{Diff}}}
\newcommand{\Cont}{\mathop{\mathrm{Cont}}}

\begin{document}


\title[Representations of diffeomorphism groups]{Irreducibility of some representations of the groups of symplectomorphisms and contactomorphisms}

\author{Łukasz Garncarek}
\address{University of Wrocław, Institute of Mathematics, pl. Grunwaldzki 2/4, 50-384 Wrocław, Poland}
\email{Lukasz.Garncarek@math.uni.wroc.pl}

\begin{abstract}
We show the irreducibility of some unitary representations of the group of symplectomorphisms and the group of contactomorphisms.
\end{abstract}

\maketitle

\section{Introduction}

Let $M$ be a smooth second-countable manifold. There exists a natural diffeo\-mor\-phism-invariant measure class on $M$, consisting of measures having positive density with respect to the Lebesgue measure in every coordinate chart. We will refer to them simply as Lebesgue measures. 

Let $\mu$ be a Lebesgue measure on $M$. For a group $G$ acting on $M$ by diffeomorphisms we may consider a series $\pmt$ of unitary representations on $L^2(M,\mu)$ given by
\begin{equation}
\pmt(\gamma)f = f\circ \gamma^{-1} \left(\frac{d\gamma_*\mu}{d\mu}\right)^{1/2+i\theta},
\end{equation} 
where $\theta\in\R$.

If a measure $\nu$ is equivalent to $\mu$, then the operator $T\colon L^2(M,\mu)\to L^2(M,\nu)$ defined by
\begin{equation}
Tf=f\left(\frac{d\mu}{d\nu}\right)^{1/2+i\theta}
\end{equation}
gives an isomorphism of representations $\pmt$ and $\Pi^\theta_\nu$. In particular, if $\mu$ is equivalent to a $G$-invariant measure, the representations $\pmt$ are equivalent for all $\theta\in\R$.

For a diffeomorphism $\phi\colon M\to M$ we define its support $\supp f$ as the closure of the set $\{p\in M : \phi(p)\ne p\}$. Compactly supported diffeo\-morphisms of $M$ form a group $\Diff_c(M)$. In~\cite{VerGelGra1982} it was proved that for an infinite measure $\mu$ the representation $\pmz$ of the group $\Diff_c(M,\mu)$ of compactly supported, measure-preserving diffeomorphisms of $M$ is irreducible. It follows that the representations $\pmt$ of the groups $\Diff_c(M,\mu)$ and $\Diff_c(M)$ are irreducible for any $\theta\in\R$. The idea of the proof is to take two functions $f,g\in L^2(M,\mu)$ and explicitly find a diffeo\-morphism $\phi$ such that $\la f, \pmz(\phi)g \ra \ne 0$, thus showing that $f$ and $g$ cannot lie in two distinct orthogonal invariant subspaces.

Representations of various subgroups of the group of diffeomorphisms are also studied in~\cite{Ismagilov}.

The purpose of this note is to present an enhancement of the argument from~\cite{VerGelGra1982}, and apply it to classical groups of diffeomorphisms: the group of symplectomorphisms and the group of contactomorphisms.

\section{Convolution on the Heisenberg group }

On $\R^n$ the following theorem holds (see Theorem 4.3.3 in~\cite{Hormander} for a proof of a more general result):
\begin{thm}\label{thm:Rn_conv}
If $f,g\in L^1(\R^n)$ are compactly supported and nonzero, then $f*g$ is nonzero.
\end{thm}
\begin{proof}
Let $\hat{h}(\xi)=\int h(x)e^{-ix\xi}\,dx$ denote the Fourier transform of $h\in L^1(\R^n)$. Suppose that $f*g=0$. As $f$ and $g$ are compactly supported, their Fourier transforms extend to entire functions. Since $\hat f \hat g = \widehat{f*g} = 0$ on $\R^n$, it follows by holomorphicity that $\hat f \hat g = 0$ on $\C^n$, and either $\hat f$ or $\hat g$ must vanish. This contradicts the assumption that $f$ and $g$ are nonzero.
\end{proof}

In this section we will prove an analogue of this theorem for square-integrable functions on the Heisenberg group.

\subsection{The Heisenberg group}
Let $n$ be a positive integer. The multiplicative group of all matrices of the form
\begin{equation}
\begin{pmatrix}
1 & \seq{x}^T & z \\
0 & I_n & \seq{y} \\
0 & 0 & 1
\end{pmatrix},
\end{equation}
where $z\in \R$, $\seq{x},\seq{y}\in\R^n$, and $I_n$ denotes the $n\times n$ identity matrix, is called the Heisenberg group $H_n$. It is a unimodular Lie group diffeomorphic with $\R^{2n+1}$, and its Haar measure is the $(2n+1)$-dimensional Lebesgue measure. We will identify $H_n$ with $\R^{2n+1}$ as manifolds. The convolution of functions $f,g\in L^1(H_n)$ will be denoted $f *_H g$.

\subsection{Convolution of compactly supported functions on $H_n$}
Let $f\in L^1(\R)$. Define
\begin{equation}
Tf(x)=\int_{-\infty}^x f(t)\,dt.
\end{equation}
If $f\in L^2(\R)$ is supported in $[a,b]$, then it is integrable; furthermore, if $\int f(t)\,dt = 0$, then $\supp Tf \subseteq [a,b]$ and we may write
\begin{equation}\label{eq:K_restr}
Tf(x) = \int f(t)K_{[a,b]}(t,x)\,dt,
\end{equation} 
where
\begin{equation}
K_{[a,b]}(t,x) = \begin{cases}
1 \qquad\text{for $a\leq t\leq x\leq b$},\\
0 \qquad\text{otherwise}.
\end{cases}
\end{equation}
Hence, $Tf \in L^2(\R)$ and $\snorm{Tf}_2 \leq \snorm{K_{[a,b]}}_2\snorm{f}_2$, where $\snorm{K_{[a,b]}}_2$ stands for the $L^2$-norm of $K_{[a,b]}\in L^2(\R^2)$. We may iterate the process of applying $T$ to $f$ as long as it yields a function integrating to $0$. The next lemma shows that unless $f=0$, this process terminates.

\begin{lem}\label{lem:int}
If $f\in L^2(\R)$ is nonzero and compactly supported, then there exists $k\geq 0$ such that $T^kf\in L^2(\R)$ and $\int T^kf(x)\,dx \ne 0$.
\end{lem}

\begin{proof}
If there is no such $k$, then $T^kf\in L^2(\R)$ and $\int T^kf(x)\,dx = 0$ for all $k$. Suppose this is the case. We may assume that $\supp f\subseteq [0,1]$, and replace $T$ with a bounded operator of the form \eqref{eq:K_restr} with kernel $K_{[0,1]}$. 

Since $f$ is compactly supported, $\hat f$ extends to an entire function on $\C$. We now have
\begin{equation}
  \widehat{T^kf}(\xi)=(i\xi)^{-k}\hat f(\xi),
\end{equation}
and by the Plancherel theorem
\begin{equation}
  4\pi^2\snorm{T^kf}_2^2 = \snorm{\widehat{T^kf}}_2^2 \geq \int_{-1}^1\abs{\hat{f}(\xi)}^2 \,d\xi
\end{equation}
But $\snorm{T}\leq\snorm{K_{[0,1]}}_2 < 1$, so the left-hand side of the above inequality can be made arbitrarily small. Therefore $\hat f= 0$, as it is an entire function vanishing on $[-1,1]$. This contradicts the assumption that $f$ is nonzero.
\end{proof}

Let $f\in L^2(H_n)$ be compactly supported. Define $Sf\in L^1(\R^{2n})$ by
\begin{equation}
  Sf(\seq{x},\seq{y}) = \int_{\R} f(\seq{x},\seq{y}, z) \,dz.
\end{equation}
If $Sf=0$, we may also define $Tf\in L^2(H_n)$  by
\begin{equation}
  Tf(\seq{x},\seq{y},z)=\int_{-\infty}^z f(\seq{x},\seq{y},t) \,dt
\end{equation}

The proof of the next lemma consists of a straightforward application of the Fubini theorem:
\begin{lem}\label{lem:Hn_conv}
If $f,g\in L^2(H_n)$ are compactly supported, then
\begin{enumerate}
  \item $S(f*_Hg) = S{f} * S{g}$,
  \item if $S{f}=0$, then $(Tf) *_H g = T(f *_H g)$,
  \item if $S{g}=0$, then $f *_H (Tg) = T(f *_H g)$.
\end{enumerate}
\end{lem}

\begin{thm} \label{thm:Hn_conv}
  If $f,g\in L^2(H_n)$ are compactly supported and nonzero, then $f *_H g \ne 0$.
\end{thm}
\begin{proof}
By Lemma~\ref{lem:int} there exist minimal $k$ and $l$ such that $ST^kf, ST^lg\in L^1(\R^{2n})$ are nonzero and compactly supported. From Lemma~\ref{lem:Hn_conv} and Theorem~\ref{thm:Rn_conv} we obtain
\begin{equation}
  ST^{k+l}(f*_H g) = S(T^kf *_H T^lg) = ST^kf*ST^lg \ne 0,
\end{equation}
which implies that $f*_H g \ne 0$.
\end{proof}

\section{Symplectic manifolds}

\subsection{Symplectic manifolds} 

Let $M$ be a symplectic manifold, that is a $2n$-dimen\-sio\-nal manifold equipped with a nondegenerate closed 2-form. A symplectomorphism of $(M,\omega)$ is a diffeomorphism $\phi\in\Diff(M)$ satisfying $\phi^*\omega=\omega$. The group of all compactly supported symplectomorphisms will be denoted by $\Sympl_c(M,\omega)$. Since $\omega$ is nondegenerate, $\omega^n$ defines a positive measure $\mu$ on $M$, invariant under the action of $\Sympl_c(M,\omega)$. 

A standard example of a symplectic manifold is $\R^{2n}$ endowed with the symplectic form $\omega_0=\sum_{i=1}^n dx^i\wedge dy^i$. It is a theorem of Darboux that any symplectic manifold is locally symplectomorphic to $(\R^{2n},\omega_0)$:
\begin{thm}\label{thm:sympl_Darboux}
For every $p\in M$ there exists a chart $\phi\colon U\to \R^{2n}$ centered at $p$, such that $\omega\res{U} = \phi^*\omega_0$.
\end{thm}
\begin{proof}
See~\cite{DaSilva}, Theorem 8.1.
\end{proof}
The chart satisfying the conditions of Theorem~\ref{thm:sympl_Darboux} is called a Darboux chart. The pushforward of $\mu$ through a Darboux chart is the standard Lebesgue measure, up to a constant factor.

The flow $\Fl^X_t$ of a complete vector field $X\in\X(M)$ consists of symplectomorphisms if and only if 
\begin{equation}\label{eq:lie_sympl}
\Ld_X \omega = 0.
\end{equation}
There is an easy way to produce such vector fields. Namely, consider a compactly supported smooth function $f\in C^\infty(M)$. Since $\omega$ is nondegenerate, there exists a unique vector field $X_f\in\X(M)$ such that $X_f \ip \omega = df$, and it is not hard to show that this field satisfies~\eqref{eq:lie_sympl}. 

For more information on symplectic manifolds see~\cite{DaSilva} and~\cite{SalamonMcDuff}.

\subsection{The representation $\pmz$ of $\Sympl_c(M,\omega)$}
As $\mu$ is a $\Sympl_c(M,\omega)$-invariant measure, the only interesting representation is $\pmz$, taking the form
\begin{equation}
\pmz(\gamma)f = f\circ\gamma^{-1}.
\end{equation}
Notice that the space of constant square-integrable functions is $\pmz$-invariant. It is nontrivial when $\mu(M)<\infty$. Let us denote its orthogonal complement by $\H$. 

\begin{thm} \label{thm:sympl_irr}
The representation $\pmz$ of the group $\Sympl_c(M,\omega)$ on the space $\H$ is irreducible.
\end{thm}

\begin{lem} \label{lem:sympl_translat}
Let $p\in M$ and let $\phi\colon U\to\R^{2n}$ be a Darboux chart centered at $p$. Then there exist $r>0$ and for every $x\in B(0,2r)$ a symplectomorphism $\tau_x \in \Sympl_c(U,\omega\res{U}) \subseteq \Sympl_c(M,\omega)$ such that
\begin{enumerate}
\item $\Cl{B(0,3r)} \subseteq \phi[U]$,
\item $\phi\tau_x\phi^{-1}(y) = y+x$ for all $y\in B(0,r)$ .
\end{enumerate}
\end{lem}
\begin{proof}
Take $r>0$ satisfying~(1) and a bump function $h\in\Ci(\R^{2n})$ supported in $\phi[U]$ and equal to $1$ on $\Cl{B(0,3r)}$. On $\R^{2n}$ there exists a linear function $f$ such that $X_f=x$ is a constant field. Then $X_{fh} = x$ on $\Cl{B(0,3r)}$ and $\supp X_{fh} \subseteq \phi[U]$. The desired symplectomorphism is $\tau_x = \phi^{-1}\Fl^{X_{fh}}_1\phi$.
\end{proof}

By using a standard argument we obtain the following well-known corollary:

\begin{cor} \label{cor:sympl_transitive}
The action of $\Sympl_c(M,\omega)$ on $M$ is $k$-transitive for all $k\geq 1$.
\end{cor}

\begin{lem} \label{lem:sympl_small_supp}
Let $\phi\colon U\to\R^{2n}$ be a Darboux chart. Then for every nontrivial $\pmz$-invariant subspace $\H_0$ of $\H$, there exists $f\in\H_0$ such that $f\ne 0$ and $\supp f\subseteq U$.
\end{lem}
\begin{proof}
We may assume that $0\in U = \phi[U]\subseteq\R^{2n}$. Let $r>0$ be as in Lemma~\ref{lem:sympl_translat}. Take a nonzero $g\in\H_0$. The $2$-transitivity of $\Sympl_c(M,\omega)$ allows us to assume without loss of generality that there exists $c\in\R$ such that the sets $A=\set{p\in B(0,r) : \Re g(p) < c}$ and $B=\set{p\in B(0,r) : \Re g(p) > c}$ both have positive measure. By the Lebesgue density theorem there exist $a\in A$ and $b\in B$ with the property that $A$ (resp. $B$) has Lebesgue density $1$ at $a$ (resp. $b$). Lemma~\ref{lem:sympl_translat} asserts the existence of a symplectomorphism $\tau=\tau_{b-a}$ that takes $a$ onto $b$ and preserves the Lebesgue density on $B(0,3r)$. The function $f=g-\pmz(\tau)g\in\H_0$ then satisfies the conclusion of the lemma.
\end{proof}

\begin{proof}[Proof of Theorem~\ref{thm:sympl_irr}] 
Suppose that $\H=\H_0\oplus\H_0^\perp$ is a nontrivial decomposition into $\pmz$-invariant subspaces. Let $\phi\colon U\to\R^{2n}$ be a Darboux chart, and let $r>0$ and $\tau_x\in\Sympl_c(M,\omega)$ be as in Lemma~\ref{lem:sympl_translat}. Without loss of generality assume that $U=\phi[U]\subseteq\R^{2n}$. By Lemma~\ref{lem:sympl_small_supp} we may choose nonzero $f\in\H_0$ and $g\in\H_0^\perp$ supported in $B(0,r)$.
We have
\begin{equation}
\la f, \pmz(\tau_x)g \ra = \int_{B(0,r)}f(y)\conj{g(\tau_x^{-1}(y))}\,dy = f*g^*(x),
\end{equation}
where $g^*(y)=\conj{g(-y)}$. But from Theorem~\ref{thm:Rn_conv} we know that this is nonzero for some $x\in\supp f*g^* \subseteq B(0,2r)$. We obtain a contradiction, since $\pmz(\tau_x)g\in\H_0^\perp$.

\end{proof}

\section{Contact manifolds}

\subsection{Contact manifolds}
Let $\dim M = 2n+1$.  A contact form on $M$ is a 1-form $\alpha\in\Omega^1(M)$ suct that $\alpha\wedge(d\alpha)^n$ is a volume form. Consider a $2n$-dimensional distribution $\xi\leq TM$. There exists an open cover $\U=\set{U_i}$ of $M$, such that for every $U\in\U$ the restriction $\xi\res{U}$ is the kernel of a 1-form $\alpha_U\in\Omega^1(U)$. If the forms $\alpha_U$ are contact forms, we call $(M,\xi)$ a contact manifold. Unless $\xi$ is the kernel of a globally defined contact form, there is no distinguished measure on $M$.

Assume for the rest of this section that $(M,\xi)$ is a contact manifold. A contactomorphism of $(M,\xi)$ is a diffeomorphism $\phi\in\Diff(M)$, such that $\phi_*\xi =\xi$. The group of compactly supported contactomorphisms will be denoted by $\Cont_c(M,\xi)$. 

An example of a contact manifold is the Heisenberg group $H_n$ with the distribution $\xi=\ker\alpha_0$, where $\alpha_0 = dz-\sum_i y^i dx^i$ is a right-invariant form on $H_n$.

There is an analogue of Darboux theorem for contact manifolds:
\begin{thm}\label{thm:cont_Darboux}
For every $p\in M$ there exists a chart $\phi\colon U\to H_n$ centered at $p$, such that $\xi\res{U} = \ker\phi^*\alpha_0$.
\end{thm}
\begin{proof}
See~\cite{Geiges}, Theorem 2.5.1.
\end{proof}

Let $U\subseteq M$ be such that $\xi\res{U}=\ker \alpha$ for some $\alpha\in\Omega^1(U)$. There exists a unique vector field $R\in\X(U)$ such that $\alpha(R) = 1$ and $R\ip d\alpha = 0$, called the Reeb vector field. If $X\in\X(U)$ is a complete vector field, then its flow $\Fl^X$ consists of contactomorphisms if and only if
\begin{equation}\label{eq:lie_cont}
\Ld_X\alpha = u\alpha
\end{equation}
for some $u\in\Ci(U)$. If we take any $f\in\Ci(U)$, by nondegeneracy of $d\alpha$ there exists $X_f\in\X(U)$ satisfying $\alpha(X_f) = f$ and $X_f\ip d\alpha = df(R)\alpha - df$. These conditions imply equality~\eqref{eq:lie_cont}. On the other hand, if $X$ satisfies~\eqref{eq:lie_cont}, then it is of the form $X_f$ for $f=\alpha(X)$.

For more information on contact manifolds see~\cite{Geiges}.

\subsection{Representations of $\Cont_c(M,\xi)$}

\begin{lem} \label{lem:cont_translat}
Let $p\in M$ and let $\phi\colon U\to H_n$ be a Darboux chart centered at $p$. Then there exist an open set $V\subseteq H_n$, a convex open neighborhood $W$ of $0$ in the Lie algebra of $H_n$, and for every $x\in \exp[W]$ a contactomorphism $\rho_x \in \Cont_c(U,\xi\res{U}) \subseteq \Cont_c(M,\xi)$ such that
\begin{enumerate}
\item $0 \in V \subseteq VV \subseteq \exp[W] \subseteq \Cl{V\exp[W]} \subseteq \phi[U]$,
\item $\phi\rho_x\phi^{-1}(y) = yx$ for all $y\in V$ .
\end{enumerate}
\end{lem}

\begin{proof}
Existence of $V$ and $W$ satisfying (1) is obvious. Let $x=\exp v$, where $v\in W$. Then $v$ extends to a left-invariand vector field $X\in\X(H_n)$, and $\Fl^X_t = R_{\exp{tv}}$, where $R_y$ is the right multiplication by $y$. If $f=h\alpha_0(X)$, where $h\res{V\exp[W]}=1$ and $\supp h\subseteq \phi[U]$, then $X_f = X$ on $V\exp[W]$. The contactomorphism $\rho_x = \phi^{-1}\Fl^{X_f}_1\phi$ satisfies condition (2).
\end{proof}

\begin{cor} \label{cor:cont_transitive}
The action of\/ $\Cont_c(M,\xi)$ on $M$ is $k$-transitive for all $k\geq 1$.
\end{cor}

\begin{lem} \label{lem:cont_small_supp}
Let $\phi\colon U\to H_n$ be a Darboux chart. Then for every nontrivial $\pmt$-invariant $\H_0 \leq L^2(M,\mu)$, there exists $f\in\H_0$ such that $f\ne 0$ and $\supp f\subseteq U$.
\end{lem}
\begin{proof}
Without loss of generality assume that $0\in U\subseteq H_n$ and $\xi\res{U} = \ker \alpha_0$. Let $\delta_t(\seq{x},\seq{y},z) = (e^t\seq{x},e^t\seq{y},e^{2t}z)$ be the flow of the field $X=(\seq{x},\seq{y},2z)$. We have $\delta_t^*\alpha_0=e^{2t}\alpha_0$, so $X=X_g$ for some function $g\in\Ci(H_n)$.

There exist $V=B(0,r)\subseteq \Cl{V} \subseteq U$ and a function $h$ supported in $U$, such that $h\res{V}=g\res{V}$. Let $\psi_t = \Fl^{X_h}_t$. Then $\psi_t\res{V} = \delta_t\res{V}$ for $t<0$. Now, by transitivity of $\Cont_c(M,\xi)$, we may take a nonzero $f\in\H_0$ such that $\supp f \cap V \ne \emptyset$. Since
\begin{equation}
\int_V \abs{\pmt(\psi_t)f}^2\,d\mu = \int_{\psi_{-t}[V]} \abs{f}^2 \,d\mu \xrightarrow[t\to\infty]{} 0,
\end{equation}
there exists $t>0$ such that $f-\pmt(\psi_t)f$ satisfies the conclusion of the lemma.
\end{proof}

Now, fix a Darboux chart $\phi\colon U\to H_n$ and a Lebesgue measure $\mu$ on $M$, such that $0\in \phi[U]$ and $\phi_*\mu$ is the standard Lebesgue measure on $\phi[U]\subseteq\R^{2n+1}$.

\begin{thm} \label{thm:cont_irr}
For every $\theta\in\R$ the representation $\pmt$ of $\Cont_c(M,\xi)$ on the space $L^2(M,\mu)$ is irreducible.
\end{thm}
\begin{proof}
The proof is analogous to the proof of Theorem~\ref{thm:sympl_irr}. Lemma~\ref{lem:cont_translat} gives us $V\subseteq U$ and contactomorphisms $\rho_x$, such that for $f$ and $g$ supported in $V$ the matrix coefficient $\la f, \pmt(\rho_x)g \ra$ is nonzero for some $\rho_x$ because of Theorem~\ref{thm:Hn_conv}.
\end{proof}

\bibliographystyle{plain}
\bibliography{reps1}

\begin{thebibliography}{1}

\bibitem{DaSilva}
Ana~Cannas da~Silva.
\newblock {\em Lectures on Symplectic Geometry}.
\newblock Springer, 2001.

\bibitem{Geiges}
Hansj{\"o}rg Geiges.
\newblock {\em An Introduction to Contact Topology}.
\newblock Cambridge University Press, 2008.

\bibitem{Hormander}
Lars H{\"o}rmander.
\newblock {\em The Analysis of Linear Partial Differential Operators I}.
\newblock Springer-Verlag, 1990.

\bibitem{Ismagilov}
R.~S. Ismagilov.
\newblock {\em Representations of Infinite-Dimensional Groups}.
\newblock American Mathematical Society, 1996.

\bibitem{SalamonMcDuff}
Dusa McDuff and Dietmar Salamon.
\newblock {\em Introduction to Symplectic Topology}.
\newblock Oxford University Press, 1998.

\bibitem{VerGelGra1982}
A.~M. Vershik, I.~M. Gel'fand, and M.~I. Graev.
\newblock Representations of the group of diffeomorphisms.
\newblock In {\em Representation theory: selected papers}. Cambridge University
  Press, 1982.

\end{thebibliography}

\end{document}